\documentclass[11pt,twoside]{article}

\usepackage{a4wide}
\usepackage{amsfonts}
\usepackage{amssymb}
\usepackage{amsmath}
\usepackage{graphicx}
\usepackage{xcolor}
\usepackage{tikz}

\newcommand{\ignore}[1]{}

\def\@begintheorem#1#2{\par\bgroup{\sc #1\ #2. }\it\ignorespaces}
\def\@opargbegintheorem#1#2#3{\par\bgroup{\sc #1\ #2\ (#3). } \it\ignorespaces}
\def\@endtheorem{\egroup}
\newtheorem{theorem}{Theorem}[section]
\newtheorem{corollary}[theorem]{Corollary}
\newtheorem{lemma}[theorem]{Lemma}

\newtheorem{example}[theorem]{Example}
\newtheorem{proposition}[theorem]{Proposition}
\newtheorem{question}[theorem]{Question}
\newtheorem{definition}[theorem]{Definition}
\newcommand{\bt}[1]{\begin{theorem}\label{#1}}
\newcommand{\bc}[1]{\begin{corollary}\label{#1}}
\newcommand{\bl}[1]{\begin{lemma}\label{#1}}
\newcommand{\be}[1]{\begin{example}\label{#1}}
\newcommand{\bp}[1]{\begin{proposition}\label{#1}}
\newcommand{\bq}[1]{\begin{question}\label{#1}}
\newcommand{\ba}[1]{\begin{algorithm}\rm\label{#1}}
\newcommand{\bd}[1]{\begin{definition}\rm\label{#1}}{\normalsize }
\newcommand{\bpr}{\noindent {\em Proof. }}
\newcommand{\et}{\end{theorem}}
\newcommand{\ec}{\end{corollary}}
\newcommand{\el}{\end{lemma}}
\newcommand{\ee}{\end{example}}
\newcommand{\ep}{\end{proposition}}
\newcommand{\eq}{\end{question}}
\newcommand{\ed}{\end{definition}}
\newcommand{\epr}{{\ \vbox{\hrule\hbox{%
\vrule height1.3ex\hskip0.8ex\vrule}\hrule}}\\\par}
\newcommand{\mepr}{{\ \ \ \vbox{\hrule\hbox{%
\vrule height1.3ex\hskip0.8ex\vrule}\hrule}}}

\def\R{\mathbb{R}}

\begin{document}

\title{\bf Quadratic Degree Sequence Optimization \\ and the\\ Critical Roots of a Graph}

\date{}

\author{Fr\'ed\'eric Meunier
\thanks{CERMICS, ENPC, Institut Polytechnique de Paris. Email: frederic.meunier@enpc.fr}
\and
Shmuel Onn
\thanks{\small Technion - Israel Institute of Technology. Email: onn@technion.ac.il}
\thanks{Corresponding author}}

\maketitle

\begin{abstract}
The {\em degree sequence optimization problem} is to find a subgraph of a
given graph which maximizes the sum over all vertices of a given function
evaluated at the subgraph degree of that vertex. Here we study this problem
and its complexity for quadratic functions. In particular, we introduce the
{\em critical roots of a graph}, and show they define intervals over which the optimal
value of the problem, as the quadratic root varies, is convex piecewise affine.\break

\end{abstract}

\section{Introduction}
\label{introduction}

Let $H=(V,F)$ be a given $n$-graph, that is, with $|V|=n$ vertices, and let
$f:[n]\rightarrow\R$ be a given function, where $[n]:=\{0,1,\dots,n-1\}$.
Here, unless otherwise specified, by a {\em subgraph} $G\subseteq H$ we mean a subgraph $G=(V,E)$ of $H$
with the same vertex set $V$ and a subset of the edges $E\subseteq F$.
We write $G=\emptyset$ for the empty subgraph $G=(V,\emptyset)$.
We denote the degree of vertex $v\in V$ in $G\subseteq H$ by $d_v(G)$.
The function $f$ induces a function on subgraphs by
$$f(G)\ :=\ \sum_{v\in V}f(d_v(G))\ .$$

We are interested in the following optimization problem denoted {\bf DSO} for brevity.

\vskip.2cm\noindent{\bf Degree Sequence Optimization.}
Given $H=(V,F)$ and $f:[n]\rightarrow\R$, solve
$${\bf DSO}:\quad \max\ \left\{f(G)\ :\ G\subseteq H\right\}\ .$$

This problem has been studied by several authors, in particular in \cite{AS},
where the function $f$ was assumed to be concave, and in \cite{DLMO},
where the graph $H$ was assumed to be complete.

\vskip.2cm
By interpolation we may and from here on do assume
that $f$ is polynomial, see Section \ref{restriction}.

\be{matching}{\bf Matching.}
Let $H$ be a given graph and $f(x):=-(x-1)^2$. Then the optimal value of DSO is zero
if and only if $H$ has a perfect matching. Indeed, $f(G)\leq0$ for any $G\subseteq H$,
with equality if and only if $d_v(G)=1$ for all $v\in V$, that is, $G$ is a perfect matching.
\ee

While DSO is trivially solvable for polynomials of degree at most $1$,
in Section \ref{restriction} we show that it is NP-hard already
for the degree 3 polynomial $f(x)=-3nx^3+18nx^2-(27n-1)x-3$.

So we restrict attention to quadratic polynomials $f(x):=a_2x^2+a_1x+a_0$.
If $a_2<0$, as in Example \ref{matching}, then $f$ is a concave function.
In this case DSO can be solved in polynomial time \cite{AS,DO}, see
Theorem \ref{minconvex} in Section \ref{convex} where we give a simplified proof of this fact.

So we assume $a_2>0$. Then, see Section \ref{restriction}, DSO with $f$ is equivalent
to DSO with a quadratic function of the form $x^2-ax=x(x-a)$, which is determined
by a single given real number $a$, its nontrivial root. So from here on we restrict attention
to the following problem.

\vskip.2cm\noindent{\bf Quadratic DSO.}
Given a graph $H=(V,F)$ and a real number $a\in\R$, solve
$${\bf QDSO}:\quad f^*(H,a)\ :=\
\max\ \left\{f(G)=\sum_{v\in V}d_v(G)\left(d_v(G)-a\right)\ :\ G\subseteq H\right\}\ .$$

\vskip.2cm
We say that connected $H$ is {\em $a$-easy} if either $G=\emptyset$ or $G=H$ is optimal
for QDSO, and arbitrary $H$ is $a$-easy if each of its connected components is. In
Section \ref{easy} we give various sufficient conditions for $H$ to be $a$-easy. In
particular, we conclude in Corollary \ref{corollary_easy} that for $a\leq 2$, every graph
$H$ is $a$-easy; for $a\leq 2.5$, every graph $H$ with no vertex of degree $2$ is $a$-easy;
and for $a\leq 3$, every graph $H$ with no vertex of degree $2$ or $3$ is $a$-easy.\break

\vskip-.2cm
Next, we make the following definition.
\bd{critical}
A {\em critical root} of graph $H$ is a real number $a$ such that the QDSO problem over $H$
and quadratic $f(x)=x(x-a)$ has optimal subgraphs $G_1=(V,E_1)$, $G_2=(V,E_2)$ with $|E_1|\neq|E_2|$.
The {\em height} of graph $H$ is the (finite) number of distinct critical roots of $H$.
\ed

For instance, nonempty $r$-regular graphs have height $1$ with $a=r$ the only critical root.

\be{regular}
Let $H\neq\emptyset$ be $r$-regular. For $a\leq r$ we have $f^*(H,a)=nr^2-nra$ with $G=H$ optimal. For
$a\geq r$ we have $f^*(H,a)=0$ with $G=\emptyset$ optimal. The only critical root is $a=r$.
\ee

In Section \ref{affinity} we study the behaviour of the optimal value $f^*(H,a)$
and optimal solutions of QDSO over a fixed $H$ as the root $a$ of the
polynomial $f(x)=x(x-a)$ varies. In particular, we show in
Theorem \ref{piecewise} that $f^*(H,a)$ is convex piecewise affine in $a$,
and a nonempty $H$ with $m$ edges has height $1\leq h\leq m$ and critical roots
$a_1<\cdots<a_h$ inducing a decomposition of the real line into $h+1$ intervals,
on each of which $f^*(H,a)$ is affine with a fixed optimal subgraph. Moreover, we provide
an algorithm for computing the critical roots and the function $f^*(H,a)$.

\vskip.2cm
The following example, which will be discussed in detail in Section \ref{affinity},
demonstrates this.
\be{affine_example}
Consider QDSO over $H=(V,F)$ drawn below, with $V=\{1,\dots,6\}$ and\break
$F=\{12,13,14,23,24,34,45,56\}$. Its height is $h=2$ and its critical roots
are $a_1=2$, $a_2=3{1\over7}$, decomposing the real line into the intervals
$C_0=(-\infty,2]$, $C_1=[2,3{1\over7}]$, $C_2=[3{1\over7},\infty)$.
For $a\in C_0$ an optimal graph is $H$, and $f^*(H,a)=48-16a$; \textcolor{red}
{for $a\in C_1$ an optimal graph is $(V,F\setminus\{56\})$}, and $f^*(H,a)=44-14a$;
and for $a\in C_2$ an optimal graph is $\emptyset$, and $f^*(H,a)=0$.\break
\ee

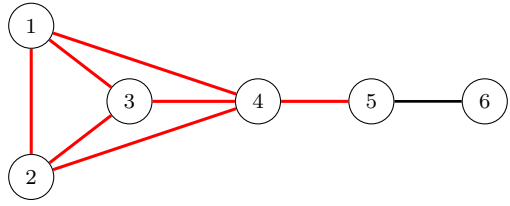
\begin{center}
\begin{tikzpicture}[every node/.style={circle,draw,font=\scriptsize},
                    % font=\small, font=\tiny,
                    % every path/.style={line width=1pt},
                    solidblack/.style={line width=1pt},
                    solidred/.style={line width=1.15pt,red},
                    dashedblack/.style={line width=1.3pt,dashed},
                    dottedblue/.style={line width=1.3pt,dotted,blue}]
  % vertices
  \node (1) at (0,1)  {1};
  \node (2) at (0,-1) {2};
  \node (3) at (1.3,0){3};
  \node (4) at (3,0)  {4};
  \node (5) at (4.5,0){5};
  \node (6) at (6,0)  {6};
  % edges
  % \draw (1) -- (2); % (if every path/.style is used)
  \draw[solidred]   (1) -- (2);
  \draw[solidred]   (1) -- (3);
  \draw[solidred]   (1) -- (4);
  \draw[solidred]   (2) -- (3);
  \draw[solidred]   (2) -- (4);
  \draw[solidred]   (3) -- (4);
  \draw[solidred]   (4) -- (5);
  \draw[solidblack] (5) -- (6);
\end{tikzpicture}
\end{center}

The proof of Theorem \ref{piecewise} uses a variant of QDSO
where the number of edges of the sought subgraph is prescribed and we just sum
degree squares. This variant was shown by Apollonio and Seb\H{o} to be NP-hard \cite{AS}.
In Section \ref{hard} we give a proof of this using a slightly different construction,
see Theorem \ref{prescribed}. Perhaps this reduction could be adapted in the future to try
and show the hardness of QDSO with unrestricted number of edges and high values of $a$.

\vskip.2cm
In Section \ref{convex} we consider DSO with concave functions.
It is known that DSO is then efficiently solvable \cite{AS,DO},
but here we give a simplified proof of this fact, see Theorem \ref{minconvex}.

\vskip.2cm
Finally, in Section \ref{trees} we consider QDSO over trees. While it was recently shown
in \cite{Onn} that DSO with any function is efficiently solvable over graphs
of bounded tree-width hence over trees, here we give a much
simpler algorithm for quadratic functions, see Theorem \ref{tree_thm}.

\section{Hardness and restriction to quadratic optimization}
\label{restriction}

We begin by proving the few simple claims stated in the introduction
leading to our restriction to degree sequence optimization with
quadratic functions of the form $f(x)=x(x-a)$.

\bp{interpolation}
Any $f:[n]\rightarrow\R$ coincides on $[n]$ with a polynomial of degree at most $n-1$.
\ep
\bpr
We can rewrite the function $f$ as the following polynomial of degree at most $n-1$,
$$f(x)\ =\ \sum_{i\in[n]}{f(i)\over g_i(i)}g_i(x)\quad \mbox{where}\quad
g_i(x):=\prod_{j\in[n]\setminus\{i\}}(x-j)\ .\mepr$$

\vskip.2cm
Next, we prove that, as noted, QDSO is NP-hard already for polynomials of degree $3$.
\bt{cubic_hard}
The DSO problem over given $H$ is NP-hard already for the cubic polynomial
$$f(x)\ :=\ -3nx^3+18nx^2-(27n-1)x-3\ .$$
\et
\bpr
A {\em cubic subgraph} is a nonempty subgraph $G\subseteq H$ in which the degree of
every $v\in V$ satisfies $d_v(G)\in\{0,3\}$.
It is known to be NP-complete to decide if $H$ has a cubic subgraph.

We claim that $\max\{f(G):G\subseteq H\}>-3n$
if and only if $H$ has a cubic subgraph.
First, rewrite $f$ as $f(x)=-3n(x-3)^2x+(x-3)$. Now, $f(0)=-3$, $f(3)=0$,
and it is easy to see that $f(x)<-3n$ for any other positive integer.
Consider any $G\subseteq H$. If $d_v(G)\notin\{0,3\}$ for some $v\in V$ then $f(G)<-3n$.
Also, if $d_v(G)=0$ for all $v\in V$ then $f(G)=-3n$. But if $d_v(G)\in\{0,3\}$
for all $v\in V$ and $d_v(G)=3$ for at least one $v\in V$ then $f(G)>-3n$.
\epr

In contrast, we show that, as noted, QDSO is trivial for polynomials of degree at most $1$.

\bp{deg_one}
DSO can be solved in polynomial time for polynomials of degree at most $1$.
\ep
\bpr
Let $f=a_1x+a_0$ be any polynomial of degree at most $1$. Then for any $G=(V,E)\subseteq H$,
$$f(G)\ =\ \sum_{v\in V}(a_1d_v(G)+a_0)\ =\ a_1\sum_{v\in V}d_v(G)+na_0\ =\ 2a_1|E|+na_0\ .$$
This implies that if $a_1\geq0$ then $G=H$ is optimal,
whereas if $a_1\leq0$ then $G=\emptyset$ is optimal.
\epr

In view of Theorem \ref{cubic_hard} and Proposition \ref{deg_one}, we restrict attention
to quadratic polynomials $f(x):=a_2x^2+a_1x+a_0$. If $a_2<0$, as in the matching
Example \ref{matching}, then $f$ is a concave function.  In this case it is known that
DSO can be solved in polynomial time, see \cite{AS,DO}. In Section \ref{convex} we give
a simplified proof of this fact, see Theorem \ref{minconvex}. So we assume $a_2>0$.
\bp{standard_quadratic}
The DSO problem with any quadratic function $f(x)=a_2x^2+a_1x+a_0$, where $a_2>0$,
is equivalent to the DSO problem with $g(x):=x(x-a)$, where $a:=-{a_1\over a_2}$.
\ep
\bpr
Let $g(x):=x(x-a)$ with $a:=-{a_1\over a_2}$.
The claim follows since for any $G\subseteq H$ we have
$$f(G)=\sum_{v\in V}\left(a_2d^2_v(G)+a_1d_v(G)+a_0\right)
=a_2\sum_{v\in V}\left(d^2_v(G)+{a_1\over a_2}d_v(G)\right)+na_0=a_2g(G)+na_0\ .\mepr$$

So from now on we may and do assume that $f(x)=x^2-ax=x(x-a)$ and is determined by a
single real number $a$, its nontrivial root, and we restrict attention to the QDSO problem.

\section{Easy graphs}
\label{easy}

In this section we show that under various assumptions on the graph $H$ and the root $a$,
QDSO with $f(x)=x(x-a)$ admits a simple optimal solution.
Here is the first such situation. Throughout we let $\delta(H)$ and $\Delta(H)$
denote the minimum and maximum vertex degrees in $H$.

\bp{observation1}
If $a\leq\delta(H)$ then $G=H$ is optimal.
If $a\geq\Delta(H)$ then $G=\emptyset$ is optimal.
If $H$ is $r$-regular, namely, $d_v(H)=r$ for all $v\in V$,
then either $G=H$ or $G=\emptyset$ is optimal.
\ep
\bpr
If $a\leq\delta(H)$ then $d_v(H)-a\geq0$ for every $v\in V$
and hence for every $G\subseteq H$ we have
$$f(G)\ =\ \sum_{v\in V}d_v(G)(d_v(G)-a)\ \leq\
\sum_{v\in V}d_v(G)(d_v(H)-a)\ \leq\ \sum_{v\in V}d_v(H)(d_v(H)-a)\ =\ f(H)\ .$$
If $a\geq\Delta(H)$ then for every $G\subseteq H$
and every $v\in V$ we have $d_v(G)-a\leq0$ and hence
$$f(G)\ =\ \sum_{v\in V}d_v(G)(d_v(G)-a)\ \leq\ 0\ =\ f(\emptyset)\ .$$
Finally, if $H$ is $r$-regular then
every $a$ satisfies either $a\leq r=\delta(H)$ or $a\geq r=\Delta(H)$.
\epr

The graphs discussed in Proposition \ref{observation1}
are easy in the sense of the following definition.

\bd{nice_def}
A connected graph $H$ is {\em $a$-easy}, where $a\in\R$, if either $G=H$ or
$G=\emptyset$ is optimal for QDSO. An arbitrary graph $H$ is {\em $a$-easy}
if each of its connected components is.
\ed
If $H$ is $a$-easy then an optimal solution $G\subseteq H$
is easy to find and is a union of some connected components of $H$.
More precisely, for every connected component $H_i$ of $H$, the
corresponding subgraph $G_i\subseteq H_i$ can be taken to be
$G_i:=H_i$ if $f(H_i)\geq 0$ and $G_i:=\emptyset$ if $f(H_i)\leq 0$.\break

The following proposition gives a simple useful criterion for $H$ to be $a$-easy.

\bp{nice_prop}
Consider QDSO over $H=(V,F)$ and $a\in\R$. If there is an optimal solution $G\subseteq H$
where every vertex $v\in V$ satisfies $d_v(G)=0$ or $d_v(G)=d_v(H)$, then $H$ is $a$-easy.
\ep
\bpr
It suffices to give the proof for a connected graph $H$.
Suppose $G=(V,E)\subseteq(V,F)=H$ is an optimal solution to QDSO where every vertex
$v\in V$ has $d_v(G)=0$ or $d_v(G)=d_v(H)$. Let $U:=\{u\in V:d_u(G)=0\}$.
Suppose for a contradiction that both $U$ and $V\setminus U$ are nonempty.
Since $H$ is connected there must be some edge $e=\{u,v\}\in F$ with $u\in U$
and $v\in V\setminus U$. Then $d_u(G)=0$ which implies $e\notin E$ whereas
$d_v(G)=d_v(H)$ which implies $e\in E$, a contradiction.
So either $U=\emptyset$ which implies $G=H$ or $U=V$ which implies $G=\emptyset$.
\epr
Call $G\subseteq H$ {\em maximal-optimal} if it is optimal for QDSO with maximum number of
edges. The following lemma establishes various useful properties of maximal-optimal subgraphs.
\bl{maximal-optimal}
Let $G=(V,E)$ be maximal-optimal for QDSO over $H=(V,F)$ and $a$. Then:
\begin{enumerate}
\item
For every edge $\{u,v\}\in F\setminus E$, we have $d_u(G)+d_v(G)<a-1$.
\item
For every vertex $v\in V$, either $d_v(G)<a-1$ or $d_v(G)=d_v(H)$.
\item
For every vertex $v\in V$ with $d_v(H)\geq 2(a-1)$, either $d_v(G)=0$ or $d_v(G)=d_v(H)$.
\item
For every vertex $v\in V$ with $d_v(H)\geq2a-1$, we have $d_v(G)=d_v(H)$.
\end{enumerate}
\el
\bpr

\vskip.2cm\noindent
1. Suppose to the contrary that $\{u,v\}\in F\setminus E$ but $d_u(G)+d_v(G)\geq a-1$.
Let $G':=(V,E')$ where $E':=E\uplus\{\{u,v\}\}$. Let $b:=d_u(G)$ and $c:=d_v(G)$.
Since $b+c\geq a-1$ we have that
\begin{eqnarray*}
f(G')-f(G) & = & \left(f(d_v(G'))-f(d_v(G))\right)+\left(f(d_u(G'))-f(d_u(G))\right) \\
& = & \left((b+1)((b+1)-a)-b(b-a)\right)+\left((c+1)((c+1)-a)-c(c-a)\right) \\
& = & 2\left((b+c)-(a-1)\right)\ \geq\ 0\ .
\end{eqnarray*}
So $G'$ is also optimal but has more edges than $G$, contradicting $G$ being maximal-optimal.

\vskip.2cm\noindent
2.
Now suppose to the contrary $a-1\leq d_v(G)<d_v(H)$ for some $v\in V$.
Since $d_v(G)<d_v(H)$ there is an edge $\{u,v\}\in F\setminus E$.
But then $d_u(G)+d_v(G)\geq a-1$ which is impossible by 1.

\vskip.2cm\noindent
3. and 4. Suppose to the contrary $v\in V$ has either $d_v(H)\geq2(a-1)$ and
$0<d_v(G)<d_v(H)$, or $d_v(H)\geq2a-1$ and $d_v(G)<d_v(H)$.
Let $U:=\{u\in V:\{u,v\}\in F\setminus E\}$, $d:=d_v(H)$, and $k:=d_v(G)$.
So $|U|=d-k$ and $d+k\geq 2a-1$. Let $G':=(V,E')$ with $E':=E\uplus\{\{u,v\}:u\in U\}$.
Consider any $u\in U$ and let $r:=d_u(G)$. Then $d_u(G')=d_u(G)+1=r+1$ and hence
$$f(d_u(G'))-f(d_u(G))\ =\ (r+1)((r+1)-a)-r(r-a)\ =\ 2r+1-a\ \geq\ 1-a\ .$$
We therefore obtain
\begin{eqnarray*}
f(G')-f(G) & = & \left(f(d_v(G'))-f(d_v(G))\right)+
\sum_{u\in U}\left(f(d_u(G'))-f(d_u(G))\right) \\
& \geq & d(d-a)-k(k-a)+(d-k)(1-a) \\
& =& d^2-k^2+(1-2a)(d-k)\ =\ (d-k)((d+k)-(2a-1))\ \geq\ 0\ .
\end{eqnarray*}
So $G'$ is also optimal but has more edges than $G$,
contradicting $G$ being maximal-optimal.
\epr
Immediate consequences of Lemma \ref{maximal-optimal}
are the following theorem and corollary.
\bt{deg_gap}
For any $a$, if $H$ has no vertex $v$ with
$2\leq d_v(H)<2(a-1)$, then $H$ is $a$-easy.
\et
\bpr
Let $G\subseteq H$ be a maximal-optimal solution. We claim that for any $v\in V$,
either $d_v(G)=0$ or $d_v(G)=d_v(H)$. Indeed, if  $d_v(H)\geq2(a-1)$
then this follows from Lemma \ref{maximal-optimal} part 3, whereas if $d_v(H)\leq 1$
then this follows trivially. So $H$ is $a$-easy by Proposition \ref{nice_prop}.
\epr
\bc{corollary_easy}
For $a\leq 2$, every graph $H$ is $a$-easy.
For $a\leq 2.5$, every graph $H$ with no vertex of degree $2$ is $a$-easy.
For $a\leq 3$, every graph $H$ with no vertex of degree $2$ or $3$ is $a$-easy.\break
\ec

In contrast, here is a graph with one vertex of degree $2$
that is not $a$-easy for all $2<a<3{1\over7}$.

\vskip.2cm\noindent
{\bf Example \ref{affine_example} continued.}
{\em
For any $2<a<3{1\over7}$, the only optimal solution of QDSO over $H=(V,F)$ where
$V=\{1,\dots,6\}$ and $F=\{12,13,14,23,24,34,45,56\}$ is $G=(V,F\setminus\{56\})$.
}

\vskip.3cm
We conclude this section with some more properties of optimal solutions for $a\geq 3$.

\bp{simplification}
Consider QDSO over $H=(V,F)$ and $a\geq3$. Suppose there is an edge $\{u,v\}\in F$ with
$d_u(H),d_v(H)\leq 2$. If $G=(V,E)$ is optimal then so is $G':=(V,E\setminus\{\{u,v\}\})$.
\ep
\bpr
Assume $\{u,v\}\in E$ else $G'=G$.
If $d_u(G)=1$ or $d_u(G)=2$ then, respectively,
$$f(d_u(G'))-f(d_u(G))\ =\ 0-1(1-a)\ =\ a-1\ >\ 0\ ,$$
$$f(d_u(G'))-f(d_u(G))\ =\ 1(1-a)-2(2-a)\ =\ a-3\ \geq\ 0\ .$$
A similar argument holds for $v$. Therefore $f(G')\geq f(G)$ and $G'$ is optimal as claimed.
\epr
If $H'$ is obtained from $H$ by repeatedly removing such edges, an optimal
$G'\subseteq H'\subseteq H$ is also optimal for $H$.
So we may assume that in $H$ there are no such edges, that is,
any vertex of degree $1$ or $2$ is adjacent only to vertices of degree at least $3$.
However, the next example shows that even under this assumption
it may be that $H$ is not $a$-easy, even when $H$ is a tree.

\be{tree_example}
The \textcolor{red}{only optimal solution} of QDSO with $a=3$ and the tree $H=(V,F)$ below,
with $V=\{1,\dots,9\}$ and $F=\{12,13,14,25,26,27,28,29\}$ is
\textcolor{red}{$G=(V,E)$ with $E=F\setminus\{13,14\}$}.
\ee

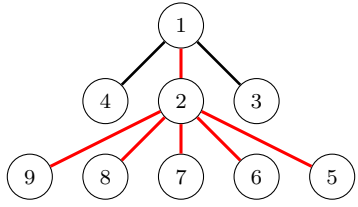
\begin{center}
\begin{tikzpicture}[every node/.style={circle,draw,font=\scriptsize},
                    solidblack/.style={line width=1pt},
                    solidred/.style={line width=1.15pt,red},
                    dashedblack/.style={line width=1.3pt,dashed},
                    dottedblue/.style={line width=1.3pt,dotted,blue}]
  % vertices
  \node (1) at (0,2)  {1};
  \node (2) at (0,1)  {2};
  \node (3) at (1,1)  {3};
  \node (4) at (-1,1) {4};
  \node (5) at (2,0)  {5};
  \node (6) at (1,0)  {6};
  \node (7) at (0,0)  {7};
  \node (8) at (-1,0) {8};
  \node (9) at (-2,0) {9};
  % edges
  \draw[solidred]  (1) -- (2);
  \draw[solidblack](1) -- (3);
  \draw[solidblack](1) -- (4);
  \draw[solidred]  (2) -- (5);
  \draw[solidred]  (2) -- (6);
  \draw[solidred]  (2) -- (7);
  \draw[solidred]  (2) -- (8);
  \draw[solidred]  (2) -- (9);
\end{tikzpicture}
\end{center}

The graph $H$ in Example \ref{tree_example} has vertices of degrees $1$ and $3$.
We now show that if $H$ has no such vertices and no adjacent vertices
of degree $2$ then $H$ is optimal for QDSO with $a=3$.

\bp{a3deg24+}
Consider QDSO over $H=(V,F)$ and $a=3$, where $\delta(H)=2$, $\Delta(H)\geq 4$,
and $H$ has no vertex of degree $3$ and no adjacent vertices of degree $2$.
Then $G=H$ is optimal.
\ep
\bpr
Let $G=(V,E)\subseteq H$ be optimal for QDSO.
Let $L:=\{\{u,v\}\in F\setminus E:d_u(H)=2\}$ and let $V_2:=\{u\in V:d_u(H)=2,\ d_u(G)=0\}$.
Let $V_k:=\{v\in V:d_v(H)=k,\ d_v(G)<k\}$ for $4\leq k\leq\Delta(H)$.
Then $2|V_2|\leq |L|\leq\sum_{k=4}^{\Delta(H)}k|V_k|$. Recall that $f(x)=x(x-3)$.
Consider any $u\in V$ with $d_u(H)=2$. If $u\in V_2$ then
$f(d_u(H))-f(d_u(G))=-2$ whereas if $u\notin V_2$ then $f(d_v(H))-f(d_v(G))=0$.
Now consider any $v\in V$ with $d_v(H)=k\geq 4$. If $v\in V_k$ then
$$f(d_v(H))-f(d_v(G))\ \geq\ k(k-3)-(k-1)((k-1)-3)\ =\ 2k-4$$
whereas if $v\notin V_k$ then $f(d_v(H))-f(d_v(G))=0$.
Then $H$ is optimal for QDSO as well, since
$$f(H)-f(G)\ \geq\ \sum_{k=4}^{\Delta(H)}(2k-4)|V_k|-2|V_2|
\ \geq\ \sum_{k=4}^{\Delta(H)}k|V_k|-2|V_2|\ \geq\ |L|-|L|=0\ .\mepr$$
We conclude this section with another useful property of maximal-optimal graphs when $a=3$.
\bp{induced}
Let $G=(V,E)$ be maximal-optimal for QDSO over $H=(V,F)$ and $a=3$.
Then $G=G[U]$ is induced by some subset $U\subseteq V$, that is, $E=\{\{u,v\}\in F:u,v\in U\}$.
\ep
\bpr
Let $U:=\{v\in V:d_v(G)\geq 1\}$. Consider any edge $\{u,v\}\in F$ with $u,v\in U$. Then we have that
$d_u(G)+d_v(G)\geq 2$ and therefore Lemma \ref{maximal-optimal} part 1 implies that $\{u,v\}\in E$.
\epr

\vskip.3cm
It would be interesting to better understand the behavior of optimal solutions
of QDSO and its complexity for small values of $a$.
By Proposition \ref{observation1} we may assume
$\delta(H)<a<\Delta(H)$. Particularly intriguing is the case of $a=3$,
and in particular the following questions. What if each degree in $H$
is $2,3,4$? What if in addition no two vertices of degree $2$ are adjacent?
\bq{a3}
What is the complexity of QDSO with $a=3$ and $d_v(H)\in\{2,3,4\}$ for all $v$?
\eq

\section{Height, critical roots, and piecewise affinity}
\label{affinity}

We now study the behaviour of the optimal value $f^*(H,a)$ and optimal solutions of QDSO
over a fixed given graph $H$ as the nontrivial root $a$ of the quadratic $f(x)=x(x-a)$ varies.

Let $H$ have $m$ edges. We need to consider variants of QDSO where the number
$0\leq e\leq m$ of edges of the sought subgraph is prescribed,
and we just {\em sum degree squares}, as follows.
$${\bf SDS}_e:\quad s_e(H)\ :=\ \max\ \left\{\sum_{v\in V}d_v^2(G) :
G=(V,E)\subseteq (V,F)=H,\ |E|=e\right\}\ , \quad e=0,\dots,m\ .$$

The following simple observation asserts that the sequence $s_0(H),\dots,s_m(H)$ is increasing.
\bp{increasing}
The optimal values of the $SDS_e$ problems satisfy $s_0(H)<\cdots<s_m(H)$.
\ep
\bpr
Suppose $G_e\subset H$ is optimal for $SDS_e$ and $e<m$. Let $G\subseteq H$ be
obtained from $G_e$ by adding any edge $uw$ which is not in $G_e$.
As the degrees of $u,w$ in $G$ are larger than in $G_e$,
$$s_{e+1}(H)\ \geq\ \sum_{v\in V}d_v^2(G)\ \geq\
\left(\sum_{v\in V}d_v^2(G_e)\right)+2\ =\ s_e(H)+2\ .\mepr$$

Fix $H$. For $e=0,\dots,m$, let $G_e\subseteq H$
be optimal for $SDS_e$, let $s_e:=s_e(H)$, and let
$$l_e\ :=\ \sup\left\{{{s_e-s_r}\over{2(e-r)}}:e<r\leq m\right\}\ \quad\mbox{and}\quad
u_e\ :=\ \inf\left\{{{s_e-s_r}\over{2(e-r)}}:0\leq r<e\right\}\ .$$
In the following theorem we show that QDSO over $H$ has an optimal solution with $e$ edges
if and only if the interval $[l_e,u_e]:=\{a\in\R: l_e\leq a\leq u_e\}$ is nonempty,
in which case $G_e$ is such an optimal solution for any $a$ in $[l_e,u_e]$.
Moreover, we show that $f^*(H,a)$ is a convex piecewise affine function in $a$ with the
intervals $[l_e,u_e]$ as its pieces. Here is the theorem, followed by a description of
an algorithm for computing the height, critical roots, $f^*(H,a)$, and optimal subgraphs,
and further followed by a demonstration via a continuation of Example \ref{affine_example}.

\bt{piecewise}
Let $H\neq\emptyset$ have $m$ edges. For some $1\leq t\leq m$,
we have reals $a_1\leq\cdots\leq a_t$, the {\em critical roots of $H$},
possibly with multiplicities, decomposing the real line into intervals,
$$C_0:=(-\infty,a_1],\ \ C_1:=[a_1,a_2],\ \ \dots,
\ \ C_{t-1}:=[a_{t-1},a_t],\ \ C_t:=[a_t,\infty),$$
integers $m=e_0>\cdots>e_t=0$, and subgraphs $H=G_{e_0},\dots,G_{e_t}=\emptyset\subseteq H$,
such that, for $i=0,\dots,t$, the subgraph $G_{e_i}$ has $e_i$ edges and is optimal
for QDSO over $H$ and any $a\in C_i$.
The optimal value $f^*(H,a)$ of QDSO is convex piecewise
affine in $a$ with the $C_i$ as its pieces.

\vskip.3cm\noindent
Moreover, letting $O_i$ be the open interval which is the
interior of $C_i$ for $i=0,\dots,t$, we have:
\begin{enumerate}
\item For $a\in C_i$: $f^*(H,a)=s_{e_i}-2e_ia$, and any SDS$_{e_i}$-optimal $G$
is also QDSO-optimal;
\item For $a\in O_i$: any QDSO-optimal $G$ has $e_i$ edges, and is also SDS$_{e_i}$-optimal.
\end{enumerate}
\et

\bpr
For $e=0,\dots,m$, let $G_e\subseteq H$ be an optimal
solution to $SDS_e$, and let $s_e:=s_e(H)$.
Consider $f(x)=x(x-a)$ with any root $a$.
For any subgraph $G\subseteq H$ with $e$ edges, we have
$$f(G)\ =\ \sum_{v\in V}\left(d_v(G)(d_v(G)-a)\right)\ =\
\sum_{v\in V}d_v^2(G)-a\sum_{v\in V}d_v(G)\ =\ \sum_{v\in V}d_v^2(G)-2ea\ .$$
Therefore
$$\max\{f(G) : G=(V,E)\subseteq (V,F)=H,\ |E|=e\}
\ =\ f(G_e)\ = \ s_e-2ea\ , \quad e=0,\dots,m\ .$$
Now,
\begin{eqnarray}
\label{max_affine}
\nonumber
f^*(H,a)
& =& \max\{f(G):G\subseteq H\} \\
& =& \max\{\max\{f(G) : G=(V,E)\subseteq (V,F)=H,\ |E|=e\}\ :\ e=0,\dots,m\} \\
\nonumber
& =& \max\{s_e-2ea\ :\ e=0,\dots,m\}\ .
\end{eqnarray}
Since $f^*(H,a)$ is the maximum of finitely many affine functions, it is indeed
convex piecewise affine. We proceed to establish the rest of the claims of the theorem.
It follows from equation \eqref{max_affine} that problem QDSO with
root $a$ has an optimal solution with $e$ edges if and only if
\begin{eqnarray}
\label{criterion1}
& & s_e-2ea\ \geq\ s_r-2ra,\ \ \ r=0,\dots,m,\ \ \ r\neq e \\
\label{criterion2}
&\Longleftrightarrow& a\leq{{s_e-s_r}\over{2(e-r)}},\ \ \ 0\leq r<e,\quad
\quad\quad a\geq{{s_e-s_r}\over{2(e-r)}},\ \ \ e<r\leq m \\
\label{criterion3}
&\Longleftrightarrow&l_e:=\sup\left\{{{s_e-s_r}\over{2(e-r)}}:e<r\leq m\right\}
\leq a\leq\inf\left\{{{s_e-s_r}\over{2(e-r)}}:0\leq r<e\right\}=:u_e
\end{eqnarray}
So QDSO has an optimal solution with $e$ edges if and only if
$[l_e,u_e]:=\{a\in\R: l_e\leq a\leq u_e\}$ is nonempty, in which case $G_e$
is such an optimal solution for any $a$ in the interval $[l_e,u_e]$.

By Proposition \ref{observation1}, $G=H$, having $m$ edges, is optimal for all sufficiently
small $a$, and $G=\emptyset$, having $0$ edges, is optimal for all sufficiently large $a$.
Thus, for some $1\leq t\leq m$ there are integers $m=e_0>\cdots>e_t=0$, which are precisely
the values of $e$ for which QDSO has an optimal solution with $e$ edges, with the corresponding
subgraphs $H=G_{e_0},\dots,G_{e_t}=\emptyset\subseteq H$, with $G_{e_i}$, the optimal
solution of $SDS_{e_i}$, also optimal for QDSO for all
$a\in C_i:=[l_{e_i},u_{e_i}]\neq\emptyset$. Note that $l_m=-\infty$, $u_0=\infty$,
and all other $l_e,u_e$ are finite, and we have in particular
$$C_0\ =\ [l_{e_0},u_{e_0}]\ =\ [l_m,u_m]\ =\
\{a\in\R\ :\ -\infty\leq a\leq u_m\}\ =\ (-\infty,u_m]\ ,$$
$$C_t\ =\ [l_{e_t},u_{e_t}]\ =\ [l_0,u_0]\ =\
\{a\in\R\ :\ l_0\leq a\leq \infty\}\ =\ [l_0,\infty)\ .$$
Since for any $a$, problem QDSO has an optimal solution with $e_i$ of edges for some $i$,
any $a$ is in some interval $C_i$, and hence the intervals $C_0,\dots,C_t$ cover $\R$.
Moreover, as claimed, for any $a\in C_i$, any SDS$_{e_i}$-optimal $G$
is also QDSO-optimal and $f^*(H,a)=s_{e_i}-2e_ia$.

On the other hand, if $a$ lies in the open interval $O_i:=(l_{e_i},u_{e_i})$, then all inequalities
in Equations \eqref{criterion1},\eqref{criterion2},\eqref{criterion3} for $e=e_i$ are
strict, so as claimed, any QDSO optimal $G$ has $e_i$ edges, and is also SDS$_{e_i}$ optimal.
This implies in particular that the $O_i$ are pairwise disjoint.

Now it remains to argue that the end points of the intervals $C_i$ are indeed ordered
as they should be. Since $f^*(H,a)$ is convex piecewise affine, its slopes are
increasing as we move over the intervals from left to right. Consider any $0\leq i<j\leq t$.
Then $e_i>e_j$, hence the slopes $-2e_i$ and $-2e_j$ of $f^*(H,a)$ over $C_i$ and $C_j$
respectively satisfy $-2e_i<-2e_j$, and therefore $C_i$ must indeed lie to the left
of $C_j$, that is, $u_{e_i}\leq l_{e_j}$. So the critical roots of $H$ are given by
$$a_1:=u_{e_0}=l_{e_1},\quad a_2:=u_{e_1}=l_{e_2},\quad
\cdots,\quad a_t:=u_{e_{t-1}}=l_{e_t}\ ,$$
some possibly appearing more than once, and hence
the {\em height of $H$} is $h:=|\{a_1,\dots,a_t\}|$.
\epr

\vskip.4cm
Critical roots may appear more than once, see Example \ref{multiple} below.
This may happen since for some $a$, there may be more than two values of $e$, say $e_i>\cdots>e_j$ for some
$0\leq i<i+1<j\leq t$, for which there are subgraphs with $e$ edges which are $a$-optimal.
Then
$$a\ =\ a_{i+1}=u_{e_i}=l_{e_{i+1}}\ =\ a_{i+2}=u_{e_{i+1}}=l_{e_{i+2}}
\ =\ \cdots\ =\ a_j=u_{e_{j-1}}=l_{e_j}\ ,$$
and the intermediate intervals are degenerate and consist of the single point $a$, that is,
$$c_{i+1}\ =\ [l_{e_{i+1}},u_{e_{i+1}}]
\ =\ \cdots\ =\ c_{j-1}\ =\ [l_{e_{j-1}},u_{e_{j-1}}]\ =\ [a,a]\ .$$
Such intermediate degenerate intervals can be removed leaving us with only $h+1$ true intervals.

\vskip.4cm
The proof of Theorem \ref{piecewise} suggests the following algorithm for computing the
height, the critical roots, the convex piecewise affine function $f^*(H,a)$, and an optimal subgraph
on each piece. Unfortunately, the algorithm is unlikely to be polynomial,
since it makes use of solving $SDS_e$, which was shown in \cite{AS} to be NP-hard,
see Section \ref{hard}. Nonetheless, it allows to compute all the needed information
in concrete situations, and is demonstrated below on Example \ref{affine_example}.

\vskip.4cm\noindent
{\bf Algorithm for the height, critical roots, $f^*(H,a)$, and optimal subgraphs.}
\begin{enumerate}

\item
Input: Nonempty graph $H=(V,F)$ with $m$ edges.

\item
For $e=0,\dots,m$ compute an optimal
solution $G_e\subseteq H$ to $SDS_e$ and let $s_e:=s_e(H)$.

\item
For $e=0,\dots,m$ compute
$$l_e\ :=\ \sup\left\{{{s_e-s_r}\over{2(e-r)}}:e<r\leq m\right\}\ ,\quad\quad
u_e\ :=\ \inf\left\{{{s_e-s_r}\over{2(e-r)}}:0\leq r<e\right\}\ .$$

\item
Let $m=e_0>\cdots>e_t=0$ be the values of $e$ for which we have $l_e\leq u_e$.

\item
Output: The critical roots $a_1:=l_{e_1},\dots,a_t:=l_{e_t}$,
height $h:=|\{a_1,\dots,a_t\}|$, intervals
$$C_0:=(-\infty,a_1],\ \ C_1:=[a_1,a_2],\ \ \dots,
\ \ C_{t-1}:=[a_{t-1},a_t],\ \ C_t:=[a_t,\infty),$$
and for $i=0,\dots,t$, the subgraph $G_{e_i}$ optimal on $C_i$,
with $f^*(H,a)=s_{e_i}-2e_ia$ on $C_i$.

\end{enumerate}

\vskip.2cm\noindent
{\bf Example \ref{affine_example} continued.}
{\em Consider again the graph $H=(V,F)$ where $V=\{1,\dots,6\}$ and $F=\{12,13,14,23,24,34,45,56\}$.
Consider the problems SDS${_e}$ for $e=0,\dots,8=m=|F|$.
Solving $SDS_e$ for $e=0,\dots,8$ by enumeration,
we find the optimal values $s_e:=s_e(H)$ to be
$$s_0=0,\ \ s_1=2,\ \ s_2=6,\ \ s_3=12,\ \ s_4=20,\ \
s_5=26,\ \ s_6=36,\ \ s_7=44,\ \ s_8=48,\ $$
with corresponding optimal subgraphs $G_e=(V,E_e)\subseteq H$
having the following sets of edges,
$$E_0=\emptyset,\ \ E_1=\{56\},\ \ E_2=\{45,56\}
,\ \ E_3=\{24,34,45\},\ \ E_4=\{14,24,34,45\}\ ,$$
$$\hskip-0.2cm E_5=\{14,23,24,34,45\},\ \ E_6=\{12,13,14,23,24,34\},\ \
E_7=\{12,13,14,23,24,34,45\},\ \ E_8=F\ .$$
For $e=0,\dots,8$ we compute the lower bound
$l_e:=\sup\left\{{{s_e-s_r}\over{2(e-r)}}:e<r\leq m\right\}$ and upper bound
$u_e:=\inf\left\{{{s_e-s_r}\over{2(e-r)}}:0\leq r<e\right\}$ on the
values of $a$, if any, for which $G_e$ is optimal, and obtain
$$l_0=3{1\over7},\ \ l_1=3{1\over2},\ \ l_2=3{4\over5},\ \ l_3=4,\ \ l_4=4,\ \
l_5=5,\ \ l_6=4,\ \ l_7=2,\ \ l_8=-\infty\ ,$$
$$u_0=\infty,\ \ u_1=1,\ \ u_2=1{1\over2},\ \ u_3=2,\ \ u_4=2{1\over2},\ \
u_5=2{3\over5},\ \ u_6=3,\ \ u_7=3{1\over7},\ \ u_8=2\ .$$
We have $l_e\leq u_e$ for $e_0=8$, $e_1=7$, $e_2=0$. The critical roots are
$a_1=l_{e_1}=l_7=2$ and $a_2=l_{e_2}=l_0=3{1\over7}$, giving the
height $h=2$, and the intervals $C_0=(-\infty,2]$, $C_1=[2,3{1\over7}]$,
$C_2=[3{1\over7},\infty)$. For $a\in C_0$ one optimal subgraph is $G_8=H$, and
$f^*(H,a)=48-16a$; for $a\in C_1$ one is $G_7$, and $f^*(H,a)=44-14a$;
and for $a\in C_2$ one is $G_0=\emptyset$, and $f^*(H,a)=0$.}

\vskip.4cm
By allowing many connected components we can design graphs $H$ that
have arbitrary critical roots with arbitrary multiplicities.
Following are two simple examples demonstrating it. The first example shows
that we can choose a critical root to have an arbitrary multiplicity.

\be{multiple}{\bf Roots with arbitrary multiplicity}.
Let $H=(V,F)$ be a perfect matching, that is, $V=\{1,\dots,2m\}$ and
$F=\{\{i,m+i\}:i=1,\dots,m\}$, and let $a=1$, so $f(x)=x(x-1)$.
For every $G\subseteq H$ and every $v\in V$ we have $d_v(G)\in\{0,1\}$
so $f(d_v(G))=0$, hence $f(G)=0$.

Therefore every subgraph is optimal, in particular $G_e$ defined to consist of the first $e$
edges of $F$ for $e=0,\dots,m$. The algorithm obtains $t=m$, $e_i=m-i$ for $i=0,\dots,m$,
and $a_1=\dots=a_m=a=1$. Therefore $a=1$ is the only critical root, appearing $m$ times, the
height is $h=1$, and the intervals are as follows, all but the last and first being degenerate,
$$C_0=(-\infty,1],\ \ C_1=[1,1],\ \ \dots,
\ \ C_{m-1}=[1,1],\ \ C_m=[1,\infty)\ .$$
Removing the degenerate middle intervals we are left with the two true intervals $(-\infty,1],[1,\infty)$.
\ee

\vskip.2cm
The second example shows that, given an arbitrary finite set
$A=\{a_1,\dots,a_k\}$ of positive integers, $a_1<\cdots<a_k$,
we can construct a graph whose set of critical roots is precisely $A$.

\be{all}{\bf Arbitrary critical roots.}
Given a set $A=\{a_1,\dots,a_k\}$ of positive integers with $a_1<\cdots<a_k$,
let $H$ be a graph with $k$ connected components $H_r$ where $H_r:=K_{a_r+1}$ for
$r=1,\dots,k$. An optimal $G\subseteq H$ is the union of $G_r\subseteq H_r$
optimal for $H_r$ for each $r$. Since $H_r$ is $a_r$-regular, by
Proposition \ref{observation1}, these optimal subgraphs are $G_r=K_{a_r+1}$
if $a_r>a$, $G_r=K_{a_r+1}$ or $G_r=\emptyset$ if $a_r=a$, and $G_r=\emptyset$ if $a_r<a$.
So if $a\in\{a_1,\dots,a_r\}$ then there are two optimal subgraphs of $H$ with numbers
of edges differing by $a+1\choose2$ and hence $a$ is a critical root, whereas if not,
then there is a unique optimal subgraph of $H$ and $a$ is not a critical root.
\ee

It would be interesting to better understand the critical roots of {\em connected} graphs,
and the complexity of computing critical roots. It seems likely that, deciding if
a given rational number $a$ is a critical root of a given graph $H$, is NP-complete,
so it would be interesting to design polynomial time algorithms for
computing critical roots of useful special graph classes.

\section{Hardness of QDSO with prescribed number of edges}
\label{hard}

The variant $SDS_e$ of QDSO with prescribed number of edges of the
sought graph and just summing degree squares was considered by
Apollonio and Seb\H{o} \cite{AS}. They show the following.

\bt{prescribed}
The problem $SDS_e$ of computing $s_e(H)$ is NP-hard, even for bipartite graphs.
\et

Their proof is quite sketchy and we encountered some difficulties in completely
making it work properly. Instead, we give here a slightly different construction,
where we connect each vertex $v$ to $n^2-d_v(H)$ rather than $\Delta(H)-d_v(H)$
new vertices of degree $1$, which makes the proof work out smoothly, and we find it useful to include it in this extended arxiv version.
Perhaps this reduction could be adapted in the future to try and show the
hardness of our QDSO problem, where the number of edges is unrestricted, perhaps for 
high values of $a$.

\vskip.2cm
\bpr
Consider the NP-complete problem of deciding, given an $n$-graph $H=(V,F)$ and $r\geq2$,
if $H$ has an $r$-clique. Given $H$ and $r$, construct a bipartite graph $B=(W=V\uplus U,C)$
as follows. For every edge $v_1v_2\in F$ introduce a new vertex $u\in U$
and edges $v_1u,v_2u\in C$. For every vertex $v\in V$ introduce $n^2-d_v(H)$
new vertices $u\in U$ and edges $vu\in C$. Thus, $d_v(B)=n^2$ for all $v\in V$
and $d_u(B)\in\{1,2\}$ for all $u\in U$. Let $e:=rn^2$. We claim that
$$\max\ \{f(G)=\sum_{w\in W}d^2_w(G)\ :\ G\subseteq B,\ |E|=e=rn^2\}\ \geq\ rn^4+rn^2+r(r-1)$$
if and only if $H$ has an $r$-clique, implying that the optimization problem is indeed NP-hard.

\vskip.2cm
For any $G\subseteq B$ and $i=1,2$ let $U_i(G):=\{u\in U:d_u(G)=i\}$.
So $U=U_1(B)\uplus U_2(B)$.

\vskip.2cm
Suppose first $K\subseteq V$ induces an $r$-clique in $H$.
Define $G=(W=V\uplus U,E)\subseteq B$ by $E:=\{vu\in C:v\in K\}$.
Then $d_v(G)=d_v(B)=n^2$ for each $v\in K$ and $d_v(G)=0$ for each $v\in V\setminus K$.
So $|E|=\sum_{v\in V} d_v(G)=rn^2=e$. Since $K$ induces a clique we have
$|U_2(G)|={r\choose2}$, and since $rn^2=|E|=\sum_{u\in U} d_u(G)=|U_1(G)|+2|U_2(G)|$
we have $|U_1(G)|=rn^2-2{r\choose2}$. So
$$f(G)\ =\ \sum_{v\in V}d^2_v(G)+\sum_{u\in U}d^2_u(G)\ =\
r\cdot(n^2)^2+(rn^2-2{r\choose2})\cdot1^2+{r\choose2}\cdot2^2\ =\ rn^4+rn^2+r(r-1)\ .$$

Now suppose $H$ has a no $r$-clique. Consider any $G=(V\uplus U,E)\subseteq B$
with $|E|=e=rn^2$. Order the vertices $V=\{v_1,\dots,v_n\}$ such that
$d_{v_1}(G)\geq\cdots\geq d_{v_n}(G)$. Let $v_k$ be the last vertex with $d_{v_k}(G)>0$.
For brevity write $d_i:=d_{v_i}(G)$ for all $i$. Since $d_i\leq n^2$ for all $i$
and $rn^2=|E|=\sum_{i=1}^nd_i=\sum_{i=1}^kd_i$ we have $k\geq r$. The largest
possible value of $\sum_{i=1}^nd^2_i$ is $rn^4$ attained when $k=r$ and $d_1=\cdots=d_r=n^2$.
The second largest possible value of $\sum_{i=1}^nd^2_i$ is $rn^4-2n^2+2$ attained when
$k=r+1$ and $d_1=\cdots=d_{r-1}=n^2$, $d_r=n^2-1$, $d_{r+1}=1$.

Suppose first $k>r$. Then $\sum_{v\in V}d^2_v(G)\leq rn^4-2n^2+2$ as explained above.
On the other hand, the largest possible value of
$\sum_{u\in U}d^2_u(G)=|U_1(G)|\cdot1^2+|U_2(G)|\cdot2^2$
with $|U_1(G)|+2|U_2(G)|=\sum_{u\in U}d_u(G)=|E|=rn^2$ is attained when $|U_2(G)|$
is as large as possible. As vertices $u\in U_2(B)$ correspond to edges in $F$,
we have $|U_2(G)|\leq|U_2(B)|=|F|\leq{n\choose2}$. So
$$\sum_{u\in U}d^2_u(G)\ =\ |U_1(G)|+4|U_2(G)|
\ \leq\ (rn^2-2{n\choose2})+4{n\choose2}\ =\ rn^2+n(n-1)$$
and
$$f(G)\ =\ \sum_{v\in V}d^2_v(G)+\sum_{u\in U}d^2_u(G)\ \leq\
(rn^4-2n^2+2)+(rn^2+n(n-1))\ <\ rn^4+rn^2+r(r-1)\ .$$

Now, suppose $k=r$. Then $d_1=\cdots=d_r=n^2$ hence $\sum_{v\in V}d^2_v(G)=rn^4$.
On the other hand, since among vertices of $V$ only $v_1,\dots,v_r$ have positive degree in $G$, distinct
vertices $u\in U_2(G)$ are connected to distinct pairs among $v_1,\dots,v_r$, and
since $H$ has no $r$-clique, at least one pair is not connected to any $u\in U_2(G)$,
hence $|U_2(G)|\leq{r\choose2}-1$. We conclude that
$$\sum_{u\in U}d^2_u(G)\ =\ |U_1(G)|+4|U_2(G)|
\ \leq\ (rn^2-2({r\choose2}-1))+4({r\choose2}-1)\ =\ rn^2+r(r-1)-2$$
and
$$f(G)\ =\ \sum_{v\in V}d^2_v(G)+\sum_{u\in U}d^2_u(G)\ \leq\
rn^4+rn^2+r(r-1)-2\ <\ rn^4+rn^2+r(r-1)\ .$$
So indeed $\max\{f(G):G\subseteq H,|E|=e\}\geq rn^4+rn^2+r(r-1)$
if and only if $H$ has an $r$-clique.
\epr

\section{Minimizing convex functions}
\label{convex}

Consider the generalization of DSO over an $n$-graph $H=(V,F)$ where for each vertex
$v\in V$ we are given a function $f_v:[n]\rightarrow\R$, and for any $G\subseteq H$
we now let $f(G):=\sum_{v\in V}f_v(d_v(G))$. Here we consider the problem in
minimization form, which is to find $G^*\subseteq H$ attaining
\begin{equation}\label{min_conv}
f(G^*)\ =\ \min\ \{f(G)\ :\ G\subseteq H\}\ .
\end{equation}
\be{deg_seq}{\bf Deciding degree sequences.}
The classical {\em degree sequence problem} is to decide if a given sequence $d=(d_v:v\in V)$
is the degree sequence $d=d(G):=(d_v(G):v\in V)$ of some graph $G=(V,E)$ over $V$.
Degree sequences have been studied by many authors, starting from their celebrated
effective characterization by Erd\H{o}s and Gallai \cite{EG}, see \cite{DLMO} and
the references therein. More generally, the problem of deciding if a given $d=(d_v:v\in V)$
is the degree sequence of some subgraph $G\subseteq H$ of a given graph $H=(V,F)$,
reduces to problem \eqref{min_conv} with a convex function $f_v(x):=(x-d_v)^2$ for each
$v\in V$, where the optimal objective function value is $f(G^*)=0$ if and only if there
exists $G^*\subseteq H$ with $d(G^*)=d$.\break The classical degree sequence problem
is obtained when $H=K_V$ is the complete graph on $V$.
\ee

As mentioned, when each $f_v$ is convex, even this more general problem can be solved in polynomial time \cite{DO}. Further, the variant of the problem where the number $0\leq e\leq m$ of edges of the sought subgraph is prescribed, is also solvable in polynomial time in \cite{AS,BH}.

Here we give a reduction of this problem to a minimum cost factor problem, which is simpler than the reduction in \cite{DO}. This reduction turns out to be quite similar to that in \cite{BH}, of which we learned after distributing an earlier version of the manuscript, but we find it useful to include it in this extended arxiv version, as it uses only simple graphs, with no multiple edges, and is tuned for the problem where the number of edges is unrestricted. 

The minimum cost factor problem is defined as follows.
Let $H=(V,F)$ be a graph and let $c:F\rightarrow\R$ be a cost function on its edges.
The {\em cost} of subgraph $G=(V,E)\subseteq H$ is
$c(G):=\sum_{e\in E}c(e)$. For each $v\in V$ let $D_v$ be an interval
(subset of consecutive integers) contained in $\{0,1,\dots,d_v(H)\}$.
A {\em $D$-factor} of $H$ is a subgraph $G\subseteq H$ such that
$d_v(G)\in D_v$ for all $v\in V$. The minimum cost factor problem
is to find a $D$-factor $G^*\subseteq H$ attaining
$$c(G^*)\ =\ \min\ \{c(G)\ :\ G\ \mbox{is a $D$-factor of}\ H\}\ .$$
The minimum cost factor problem is solvable in strongly polynomial time \cite[Theorem 35.2]{Sch}.

\vskip.2cm
We now reduce minimum convex degree sequence optimization to minimum cost factor.
Given input $H=(V,F)$ and convex function $f_v$ for each $v\in V$, we extend $H$ to a graph
${\bar H}=({\bar V},{\bar F})$ and introduce a cost function on its edges
$c:{\bar F}\rightarrow\R$ as follows. We introduce a new vertex $r$, and for each
$v\in V$ we introduce $d_v:=d_v(H)$ new vertices $v^1,\dots,v^{d_v}$. We let
${\bar V}:=V\uplus\{r\}\uplus\biguplus_{v\in V}\{v^1,\dots,v^{d_v}\}$.
We let ${\bar F}:=F\uplus\biguplus_{v\in V}\biguplus_{i=1}^{d_v}\{\{v,v^i\},\{r,v^i\}\}$.
Since $\sum_{v\in V}d_v=2m$ where $m:=|F|$,
we have $|{\bar V}|=n+2m+1$, $|{\bar F}|=5m$, and $d_r({\bar H})=2m$.

Now, for $v\in V$ and $i=1,\dots,d_v$ we let $c_v^i:=f_v(i)-f_v(i-1)$, we let
the cost of edge $\{r,v^i\}$ be $c(\{r,v^i\}):=c_v^i$, and we let
all other edges $e\in{\bar F}$ have zero cost $c(e):=0$.

We let $D_r:=\{0,1,\dots,2m\}$, and for each $v\in V$ we let $D_v:=\{d_v\}$
and for $i=1,\dots,d_v$ we let $D_{v^i}:=\{1\}$.
The minimum cost factor problem is then to find ${\bar G}^*\subseteq{\bar H}$ attaining
$$c({\bar G^*})\ =\
\min\ \{c({\bar G})\ :\ {\bar G}\ \mbox{is a $D$-factor of}\ {\bar H}\}\ .$$
The problem is feasible since ${\bar G}=({\bar V},{\bar E})$
with ${\bar E}:=F\uplus\biguplus_{v\in V}\biguplus_{i=1}^{d_v}\{\{r,v^i\}\}$ is a $D$-factor.

\vskip.2cm
We have the following theorem.
\bt{minconvex}
If ${\bar G^*}$ is a minimum cost factor then its restriction
$G$ to $H$ is optimal for the convex minimization problem.
In particular, the latter is solvable in strongly polynomial time.
\et
\bpr
Let ${\bar G^*}=({\bar V},{\bar E})$ be a minimum cost factor and let $G$ be its
restriction to $H$. Consider any $v\in V$. Since $D_v=\{d_v\}$ and $v$ is adjacent in $G$
to $d_v(G)$ edges, $v$ is adjacent in ${\bar G^*}$ to $d_v-d_v(G)$ of the vertices $v^i$.
Since each $D_{v^i}=\{1\}$ we have that the other $d_v(G)$ vertices $v^i$ are adjacent
in ${\bar G^*}$ to $r$. Since $f_v$ is convex we have
$c^1_v\leq c^2_v\leq\cdots\leq c^{d_v}_v$ and therefore
$\sum\{c_v^i:\{r,v^i\}\in{\bar E}\}\geq\sum_{i=1}^{d_v(G)}c_v^i$.
Noting that only edges $\{r,v^i\}$ have nonzero cost, we obtain
\begin{eqnarray}\label{eq1}
c({\bar G^*}) & =& \sum_{v\in V}\sum\{c_v^i:\{r,v^i\}\in{\bar E}\}
\ \geq\ \sum_{v\in V}\sum_{i=1}^{d_v(G)}(f_v(i)-f_v(i-1)) \\
\nonumber
& =& \sum_{v\in V}(f_v(d_v(G))-f_v(0))\ =\ f(G)-f(\emptyset)\ .
\end{eqnarray}
Now consider an optimal solution $G^*\subseteq H$ to the convex minimization problem.
Extend it to a factor ${\bar G}$ of ${\bar H}$ by adding, for each $v\in V$, the edges
$\{r,v^i\}$ for $1\leq i\leq d_v(G^*)$ and $\{v,v^i\}$ for $d_v(G^*)<i\leq d_v$.
It is indeed a factor, since for each $v\in V$ we have
$d_v({\bar G})=d_v(G^*)+(d_v-d_v(G^*))=d_v$ and $d_v^i=1$ for $i=1,\dots,d_v$.
The cost of this factor is
\begin{eqnarray}\label{eq2}
c({\bar G}) & =& \sum_{v\in V}\sum_{i=1}^{d_v(G^*)}c_v^i
\ =\ \sum_{v\in V}\sum_{i=1}^{d_v(G*)}(f_v(i)-f_v(i-1)) \\
\nonumber
& =& \sum_{v\in V}(f_v(d_v(G^*))-f_v(0))\ =\ f(G^*)-f(\emptyset)\ .
\end{eqnarray}
Combining equations \eqref{eq1} and \eqref{eq2} and $c({\bar G^*})\leq c({\bar G})$ we obtain
$$f(G)\ \leq\ f(\emptyset)+c({\bar G^*})\leq f(\emptyset)+c({\bar G})\ =\ f(G^*)$$
and conclude that $G$ is indeed also optimal to the convex minimization problem as claimed.
\epr

\be{triangle}
Consider the triangle $H=(V,F)$ with $V=\{u,v,w\}$ and $F=\{uv,uw,vw\}$
and convex functions $f_u(x)=f_w(x)=(x-1)^2$ and $f_v(x)=(x-2)^2$.
We extend $H$ to a graph ${\bar H}=({\bar V},{\bar F})$ with vertex set
${\bar V}:=\{u,v,w,r,u^1,u^2,v^1,v^2,w^1,w^2\}$ and edge set
$${\bar F}\ :=\ \{uv,uw,vw,uu^1,uu^2,ru^1,ru^2,vv^1,vv^2,rv^1,rv^2,ww^1,ww^2,rw^1,rw^2\}\ .$$
We let the edge costs $c$ be
$c(ru^1):=c(rw^1):=-1$, $c(ru^2):=c(rw^2):=1$, $c(rv^1):=-3$, $c(rv^2):=-1$,
and $c(e):=0$ for all other edges $e\in{\bar F}$.
Finally, we let the vertex intervals be
$$D_r:=\{0,1,\dots,6\},\ \ D_u:=D_v:=D_w:=\{2\},
\ \ D_{u^i}:=D_{v^i}:=D_{w^i}:=\{1\},\ i=1,2\ .$$
It is easy to see that taking all edges with negative cost and no edges with positive cost,
we can add (in a unique way) some edges with zero cost, obtaining a
\textcolor{red}{$D$-factor ${\bar G^*}=({\bar V},{\bar E})$, which is necessarily of
minimum cost, where ${\bar E}:=\{uv,vw,ru^1,uu^2,rv^1,rv^2,rw^1,ww^2\}$.
Its restriction $G=(V,E)$ to $H$, where \textcolor{red}{$E=\{uv,vw\}$},
is indeed optimal for the convex minimization problem.}
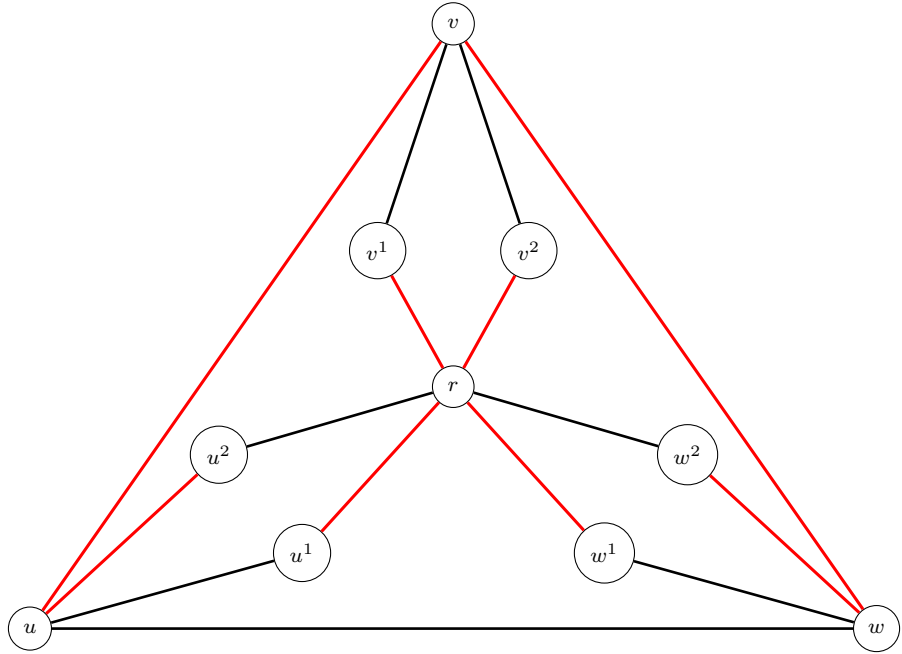
\ee

\begin{center}
\begin{tikzpicture}[every node/.style={circle,draw,font=\scriptsize},
                    solidblack/.style={line width=1pt},
                    solidred/.style={line width=1.15pt,red},
                    dashedblack/.style={line width=1.3pt,dashed},
                    dottedblue/.style={line width=1.3pt,dotted,blue}]
  % vertices
  \node (0) at (-5.6,-4)   {$u$};
  \node (1) at (-2,-3)     {$u^1$};
  \node (2) at (-3.1,-1.7) {$u^2$};
  \node (3) at (0,4)       {$v$};
  \node (4) at (-1,1)      {$v^1$};
  \node (5) at (1,1)       {$v^2$};
  \node (6) at (5.6,-4)    {$w$};
  \node (7) at (2,-3)      {$w^1$};
  \node (8) at (3.1,-1.7)  {$w^2$};
  \node (9) at (0,-0.8)    {$r$};
  % edges
  \draw[solidblack] (0) -- (1);
  \draw[solidred] (0) -- (2);
  \draw[solidblack] (3) -- (4);
  \draw[solidblack] (3) -- (5);
  \draw[solidblack] (6) -- (7);
  \draw[solidred] (6) -- (8);
  \draw[solidred] (0) -- (3);
  \draw[solidblack] (0) -- (6);
  \draw[solidred] (3) -- (6);
  \draw[solidred] (9) -- (1);
  \draw[solidblack] (9) -- (2);
  \draw[solidred] (9) -- (4);
  \draw[solidred] (9) -- (5);
  \draw[solidred] (9) -- (7);
  \draw[solidblack] (9) -- (8);
\end{tikzpicture}
\end{center}

Can we in fact also {\em maximize} arbitrary convex functions in polynomial time,
that is, find $G^*$ with $f(G^*)=\max\{f(G):G\subseteq H\}$ where
$f(G):=\sum_{v\in V}f_v(d_v(G))$ with each $f_v$ convex? This is not excluded by
Theorem \ref{cubic_hard} since the cubic polynomial therein is not convex.

\section{Trees}
\label{trees}

We now consider QDSO over trees. While it was recently shown in \cite{Onn} that DSO
with any function is efficiently solvable over graphs of bounded tree-width and
hence over trees, here we give a much simpler algorithm for quadratics.
By Corollary \ref {corollary_easy} we may assume $a>2$.

\vskip.2cm
We consider rooted trees. This does not affect DSO.
Let $T=(V,F)$ be a tree with root $r\in V$.
For each $v\in V$ let $D_v$ be the set of descendants of $v$ (including $v$)
and let $C_v\subset D_v$ be the set of children of $v$ (excluding $v$).
Let $T_v=(D_v,F_v)$ be the subtree of $T$ rooted at $v$.

Consider a tree $T$ with root $r$ and quadratic function $f(x)=x(x-a)$ with $a>2$.
For any subgraph $G\subseteq T$ we have as usual $f(G)=\sum_{v\in V}f(d_v(G))$,
and we also define $g(G):=f(d_r(G)+1)+\sum_{v\in V\setminus\{r\}}f(d_v(G))$.
We will solve recursively, over subtrees $T_v$ of $T$ from the leaves to the root,
the QDSO problem as well the following auxiliary problem {\bf AUX}.

\vskip.2cm\noindent{\bf Auxiliary problem.}
Given rooted $T$ and $f(x)=x(x-a)$ with $a>2$, solve
$${\bf AUX}:\quad \max\left\{g(G):G\subseteq T\right\}\ .$$

\vskip.2cm
Roughly speaking, for the subtree $T_v$ of $v\in V$, problem QDSO over $T_v$
gives the contribution of the vertices of $T_v$ to the value $f(G)$ of an optimal
solution $G$ of QDSO over $T$ assuming the edge from $v$ to its
parent is not present in $G$, whereas problem AUX over $T_v$ gives the contribution
assuming the edge from $v$ to its parent is present in $G$.

\vskip.2cm
As for QDSO, a {\em maximal-optimal} solution of AUX is an optimal
solution $G^+\subseteq T$ with a maximum possible number of edges.
By Lemma \ref{maximal-optimal}, for any maximal-optimal $G^*\subseteq T$ of QDSO
and any $v$, either $d_v(G^*)<a-1$ or $d_v(G^*)=d_v(T)$.
It is easy to see that a similar argument to that of Lemma \ref{maximal-optimal}
shows that the same holds for AUX, namely, for any maximal-optimal $G^+\subseteq T$
of AUX and any $v$, either $d_v(G^+)<a-1$ or $d_v(G^+)=d_v(T)$.

The following algorithm solves QDSO and AUX recursively on subtrees from the leaves
down to the root. For each $v\in V$ it determines optimal solutions
$G^*_v,G^+_v\subseteq T_v$ to the problems QDSO and AUX over the subtree $T_v$ of $T$
rooted at $v$. In what follows, the {\em union} of two graphs
$G_i=(V_i,E_i)$, $i=1,2$, is the graph $G_1\cup G_2:=(V_1\cup V_2,E_1\cup E_2)$.

\vskip.4cm\noindent
{\bf Algorithm for QDSO and AUX over trees.}
\begin{enumerate}
\item
Input: A tree $T=(V,F)$ rooted at $r\in V$, and a real number $a>2$.
\item
For each $v\in C_r$, apply the algorithm recursively to the subtree $T_v=(D_v,F_v)$ of $T$
rooted at $v$, and the real number $a$, and obtain $G^*_v$, $G^+_v$ and $f(G^*_v)$, $g(G^+_v)$.
\item
Find
$$C^*\ \in\ \arg\max\{f(C):C\subseteq C_r,\ |C|<a-1\ \mbox{or}\ C=C_r\}\ ,$$
$$C^+\ \in \ \arg\max\{g(C):C\subseteq C_r,\ |C|<a-1\ \mbox{or}\ C=C_r\}\ ,$$
where for each $C\subseteq C_r$,
$$f(C)\ :=\ |C|\left(|C|-a\right)+\sum_{v\in C}g(G^+_v)+\sum_{v\in C_r\setminus C}f(G^*_v)
\ ,\quad g(C)\ :=\ f(C)+2|C|+1-a\ ,$$
$$G(C)\ :=\ \left(\{r\}\cup C,\{\{r,v\}:v\in C\}\right)\cup
\bigcup\{G^+_v:v\in C\}\cup\bigcup\{G^*_v:v\in C_r\setminus C\}\ .$$
\item
Output: $G^*_r:=G(C^*)$, $G^+_r:=G(C^+)$ and $f(G^*_r)$, $g(G^+_r)$.
\end{enumerate}

\vskip.2cm
A total ordering of $V$ is {\em compatible} with the rooted tree $T$
if every $v\in V$ is last in the set $D_v$ of descendants of $v$ in $T$.
Note that the algorithm eventually solves recursively QDSO and AUX over all subtrees
$T_v$ in some compatible order of $V$. So when implementing the algorithm it is
convenient to choose ahead some compatible order and solve QDSO and AUX over the subtrees
$T_v$ according to this order. Then, on the recursive execution of the algorithm
over any subtree $T_v$, step 2 becomes redundant, since each $u\in C_v$ precedes $v$
in the order so $G^*_u$, $G^+_u$, $f(G^*_u)$, $g(G^+_u)$ are already available.
Choosing a compatible order can be done easily in a greedy way:
first put the leaves in any order, peel them off, next put the new resulting leaves
in any order, and so on, till at last the root is reached and added last to the order.

\vskip.2cm
Before showing that this algorithm solves QDSO and AUX in polynomial time we give an
example illustrating the algorithm with a pre choice of compatible order implemented.

\vskip.2cm\noindent
{\bf Example \ref{tree_example} continued.}
{\em
Consider $a=3$ and the tree $T=(V,F)$ with $V=\{1,\dots,9\}$, rooted at $r=1$,
$F=\{12,13,14,25,26,27,28,29\}$, and set of leaves $L=\{3,4,5,6,7,8,9\}$.
The algorithm solves recursively QDSO and AUX over all subtrees $T_v$ from the leaves
down to the root in the compatible order $3,\dots,9,2,1$, starting with the leaves,
then $2$, ending with the root $1$. We proceed to solve QDSO and AUX over all
subtrees $T_v$ in this order, obtaining $G^*_v$, $G^+_v$, $f(G^*_v)$, $g(G^+_v)$,
so these are available when considering the next vertex in order.

\begin{itemize}

\item
Leaves: at each leaf $l\in L$ we have
$G^*_l=G^+_l=(\{l\},\emptyset)$ and $f(G^*_l)=0$, $g(G^+_l)=-2$.

\item
Vertex $2$, having set of children $C_2=\{5,6,7,8,9\}\subset L$:

\begin{itemize}
\item
If $C=\emptyset$ then $f(C)=\sum_{l\in C_2}f(G^*_l)=0$ and $g(C)=f(C)-2=-2$.
\item
If $C=\{l\}$ for any $l\in C_2$ then
$$f(C)=-2+g(G^+_l)+\sum_{u\in C_2\setminus\{l\}}f(G^*_u)=-4\ ,\quad g(C)=f(C)=-4\ .$$
\item
If $C=C_2$ then $f(C)=10+\sum_{l\in C_2}g(G^+_l)=0$ and $g(C)=f(C)+8=8$.
\end{itemize}

So $C^*:=\emptyset$ and $C^+:=C_2$ are maximizers, giving
$f(G^*_2)=0$ and $g(G^+_2)=8$ where
$$G^*_2\ :=\ (\{2\},\emptyset)\cup\bigcup_{l=5}^9 G^*_l\ =\ (\{2,5,6,7,8,9\},\emptyset)\ ,$$
$$G^+_2:=\left(\{2\}\cup C_2,\{2u:u\in C_2\}\right)\cup\bigcup_{l=5}^9 G^+_l
=(\{2,5,6,7,8,9\},\{25,26,27,28,29\})\ .$$

\item
Root $1$, having set of children $C_1=\{2,3,4\}$:

\begin{itemize}
\item
If $C=\emptyset$ then $f(C)=\sum_{u\in C_1}f(G^*_u)=0$ and $g(C)=f(C)-2=-2$.
\item
If $C=\{2\}$ then $f(C)=-2+g(G^+_2)+g(G^*_3)+g(G^*_4)=6$ and $g(C)=f(C)=6$.
\item
If $C=\{3\}$ or $\{4\}$ then $f(C)=-2+g(G^*_2)+g(G^+_3)+g(G^*_4)=-4$ and $g(C)=-4$.
\item
If $C=C_1$ then $f(C)=\sum_{u\in C_1}g(G^+_u)=8-2-2=4$ and $g(C)=f(C)+4=8$.
\end{itemize}

So $C^*:=\{2\}$ are $C^+:=C_1$ are maximizers, giving $f(G^*_1)=6$, $g(G^+_1)=8$, where
$$G^*_1\ :=\ (\{1,2\},\{12\})\cup G^+_2\cup G^*_3\cup G^*_4
\ =\ (V,\{12,25,26,27,28,29\})\ ,$$
$$G^+_1\ :=\ (\{1,2,3,4\},\{12,13,14\})\cup\bigcup_{u=2}^4 G^+_u\ =\ (V,F)\ =\ T\ .$$

\end{itemize}
So the algorithm obtains
$G^*_1=(V,\{12, 25,26,27,28,29\})$, $G^+_1=T$, $f(G^*_1)=6$, $f(G^+_1)=8$.
}

\vskip.3cm
We now show that the algorithm solves QDSO and AUX in polynomial time.
\bt{tree_thm}
The algorithm above solves QDSO and AUX over trees in polynomial time.
\et
\bpr
We begin by showing that the algorithm correctly solves problems QDSO and AUX.

We use induction on the number of vertices of $T$.
If $|V|=1$ then $C^*_r=C^+_r=C_r=\emptyset$ and
$G^*_r=G^+_r=(\{r\},\emptyset)=T$ are trivially optimal with $f(G^*_r)=0$, $g(G^+_r)=1-a$.

Now assume $|V|>1$. Let $G^*_v, G^+_v$ for $v\in C_r\cup\{r\}$
be the graphs obtained by the algorithm.

Let $G=(V,E)$ be a maximal-optimal solution to QDSO. Let $C:=\{v\in C_r:rv\in E\}$.
By Lemma \ref{maximal-optimal} we have $|C|<a-1$ or $C=C_r$.
Consider any $v\in C_r$. Let $G_v:=(D_v,F_v\cap E)\subseteq T_v$
be the restriction of $G$ to $T_v=(D_v,F_v)$.
If $v\notin C$ then
$d_u(G)=d_u(G_v)$ for each $u\in D_v$ and hence $\sum_{u\in D_v}f(d_u(G))=f(G_v)$.
If $v\in C$ then $d_v(G)=d_v(G_v)+1$ and $d_u(G)=d_u(G_v)$ for each
$u\in D_v\setminus\{v\}$ and hence $\sum_{u\in D_v}f(d_u(G))=g(G_v)$.
By induction, $G^*_v$, $G^+_v$ are optimal for QDSO, AUX over $T_v$,
and hence $f(G_v)\leq f(G^*_v)$, $g(G_v)\leq g(G^+_v)$. We conclude that
\begin{eqnarray*}
f(G)&=&\sum_{v\in V}f(d_v(G))=
f(d_r(G))+\sum_{v\in C}\sum_{u\in D_v}f(d_u(G))
+\sum_{v\in C_r\setminus C}\sum_{u\in D_v}f(d_u(G)) \\
&=&|C|\left(|C|-a\right)+\sum_{v\in C}g(G_v)+\sum_{v\in C_r\setminus C}f(G_v) \\
&\leq & |C|\left(|C|-a\right)+\sum_{v\in C}g(G^+_v)+\sum_{v\in C_r\setminus C}f(G^*_v) \\
&=& f(C)\leq f(C^*) =
|C^*|\left(|C^*|-a\right)+\sum_{v\in C^*}g(G^+_v)+\sum_{v\in C_r\setminus C^*}f(G^*_v) \\
&=& f(d_r(G^*_r))+\sum_{v\in C^*}\sum_{u\in D_v}f(d_u(G^*_r))
+\sum_{v\in C_r\setminus C^*}\sum_{u\in D_v}f(d_u(G^*_r))\\
&=&\sum_{v\in V}f(d_v(G^*_r))=f(G^*_r).
\end{eqnarray*}
So $G^*_r$ output by the algorithm is also optimal for QDSO.
Now let $G=(V,E)$ be a maximal-optimal solution to AUX. Let $C:=\{v\in C_r:rv\in E\}$.
Let $G_v$ be the restriction of $G$ to $T_v$. A similar calculation below
shows that $G^+_r$ output by the algorithm is also optimal for AUX,
\begin{eqnarray*}
g(G)&=& f(d_r(G)+1)+\sum_{v\in C}\sum_{u\in D_v}f(d_u(G))
+\sum_{v\in C_r\setminus C}\sum_{u\in D_v}f(d_u(G)) \\
&=&(|C|+1)\left((|C|+1)-a\right)+\sum_{v\in C}g(G_v)+\sum_{v\in C_r\setminus C}f(G_v) \\
&\leq & (|C|+1)\left((|C|+1)-a\right)+\sum_{v\in C}g(G^+_v)+
\sum_{v\in C_r\setminus C}f(G^*_v)= f(C)+2|C|+1-a \\
&=& g(C)\leq g(C^+)= f(C^+)+2|C^+|+1-a \\
&=& (|C^+|+1)\left((|C^+|+1)-a\right)+\sum_{v\in C^+}g(G^+_v)+
\sum_{v\in C_r\setminus C^+}f(G^*_v) \\
&=& f(d_r(G^+_r)+1)+\sum_{v\in C^+}\sum_{u\in D_v}f(d_u(G^+_r))
+\sum_{v\in C_r\setminus C^+}\sum_{u\in D_v}f(d_u(G^+_r)) \\
&=& \sum_{v\in V}f(d_v(G^+_r))=g(G^+_r).
\end{eqnarray*}
The induction step follows and the proof of the correctness of the algorithm is completed.

\vskip.2cm
We continue to show that the running time of the algorithm is polynomial.

To obtain $C^*,C^+$ in step 3, we do not need to compute $f(C)$, $g(C)$ for all
relevant $C\subseteq C_r$; instead, for each integer $0<k<a-1$, we let $B_k\subseteq C_r$
consist of the $k$ elements $v\in C_r$ with largest values $g(G^+_v)-f(G^*_v)$.
Then for every $C\subseteq C_r$ with $|C|=k=|B_k|$ we have
\begin{eqnarray*}
f(C)&=&|C|\left(|C|-a\right)+\sum_{v\in C}g(G^+_v)+\sum_{v\in C_r\setminus C}f(G^*_v) \\
&=& k(k-a)+\sum_{v\in C_r}f(G^*_v)+\sum_{v\in C}\left(g(G^+_v)-f(G^*_v)\right) \\
&\leq & k(k-a)+\sum_{v\in C_r}f(G^*_v)+\sum_{v\in B_k}\left(g(G^+_v)-f(G^*_v)\right) \\
&=&|B_k|\left(|B_k|-a\right)+\sum_{v\in B_k}g(G^+_v)+
\sum_{v\in C_r\setminus B_k}f(G^*_v)\ =\ f(B_k)\ .
\end{eqnarray*}
So $B_k$ gives the best value $f(C)$ among $C\subseteq C_r$ with $|C_r|=k$.
Since for such sets we have $g(C)=f(C)+2k+1-a$, $B_k$ also gives best value $g(C)$.
So $C^*$, $C^+$ can be taken to be the ones giving best values $f(C)$, $g(C)$
among $C=\emptyset$, $C=B_k$ for $0<k<a-1$, and $C=C_r$.

It is easy to see that the time taken by step 3 is bounded by $can$ for some constant $c$.
As we just want to argue that the running time $t(n)$ is polynomial, we just show it is
bounded by $can^2$. Letting $d:=d_r(T)$ be the degree of the root,
$C_r=\{v_1,\dots,v_d\}$ its set of children, and $n_i:=|D_{v_i}|$ the number
of vertices in the subtree $T_{v_i}$ rooted at $v_i$ for $i=1,\dots,d$, we have
$$\hskip-.1cm
t(n)\ \leq\ can+\sum_{i=1}^d t(n_i)\ \leq\ ca(n+\sum_{i=1}^dn_i^2)
\ \leq\ ca(n+(\sum_{i=1}^dn_i)^2)\ =\ ca(n+(n-1)^2)\ \leq \ can^2\ .\mepr$$

\section*{Acknowledgments}

Shmuel Onn was partially supported by the Dresner chair at the Technion.

\end{document}